\magnification=\magstep1
\parindent=0mm
\parskip=0mm
\baselineskip=18pt 
\bigskipamount=18pt plus 6 pt minus 6 pt
\medskipamount=14pt plus 3 pt minus 3 pt
\smallskipamount=9pt plus 2 pt minus 2 pt
\thinmuskip=5mu
\medmuskip=6mu plus 3mu minus 2mu
\thickmuskip=6mu plus 3mu minus 2mu
\def\7{{\hskip-4pt}}
\overfullrule=0pt

\font\xivrm=cmr17

\def\R{{I\mkern-5mu R}}
\def\N{{I\mkern-5mu N}}

\def\Any{\,\cdot\,}

\def\giantrightarrow{\mathrel{\hbox to 1.5cm{\rightarrowfill}}}
\def\ghookrightarrow{\mathrel{\hbox to 1.5cm{$\lhook\mkern-9mu$\rightarrowfill}}}
\def\square{\mathord{\vbox{\hrule\hbox{\vrule\hskip 9pt\vrule height 9pt}\hrule}}}

\def\dots{{.\hskip-.05em .\hskip-.05em .}}
\frenchspacing

{\xivrm
\centerline{On the integrability in finite terms of first-order}
\centerline{differential equations according to Maximovi\v c.}
}

\bigskip
\centerline{Karl Michael Schmidt}
\centerline{Mathematisches Institut der Universit\"at,}
\centerline{Theresienstr. 39, D-80333 M\"unchen, Germany}

\bigskip
\centerline{1991 Mathematics Subject Classification: 34A46, 34A05, 01A70}

\bigskip
{\bf Abstract.}
\qquad
{\it
We give a definition of integration by quadratures of first-order ordinary
differential equations,
and recover a little known result by Maximovi\v c
which states that
a first-order ordinary differential equation can be integrated by quadratures
only if it arises from the linear equation by a diffeomorphic transformation
of the dependent variable.
In the appendix this result is applied to the linear second-order equation
in a restricted setting.
A brief outline of Maximovi\v c's life is also included.
\/}

\bigskip
{\bf Introduction.}

\bigskip
It is the fundamental reason for the continuing interest in the subject,
as well as a source of exasperation and disappointment to the beginning student,
that there is no general method or algorithm to solve a given ordinary
differential equation, in fact that there are only very few 
differential equations which can be integrated in finite terms, i.e.
whose solutions can be explicitly given in terms of elementary functions
and processes.
Several attempts at a systematic treatment of the question which equations
can be integrated in finite terms have been undertaken, beginning with
Liouville's theory (cf. [7]);
later Lie's study of transformation groups was motivated by a classification
of (partial) differential equations and integration methods for their solution
(see [5] for an account).
His ideas have subsequently given rise to the voluminous theory of
Picard-Vessiot extensions in differential algebra (cf. [6]).

In his 1885 treatise [4], V P Maximovi\v c has developed a different approach.
Based on the well-known explicit solution formula for the linear first order
equation
$$\eqalignno{
& y' + p\, y = q, & (1)
\cr
& y(x) = e^{-\int_{x_0}^x p}
         \Big(y(x_0) + \int_{x_0}^x q(t)\, e^{\int_{x_0}^t p}\, dt\Big),
  & (2)
\cr}$$
which contains the two {\it quadratures\/} $s = \int p$ and
$S = \int q\,e^{-s}$,
he considers the question under which circumstances a first-order ordinary
differential equation can be {\it integrated by quadratures\/},
i.e. its general solution can be expressed as an (elementary) function of
finitely many arbitrarily nested quadratures, their integration constants
serving as parameters for the solution manifold.
Note that in the above example the same solution formula is valid for all
(reasonable) coefficient functions $p$ and $q$: similarly, Maximovi\v c
studies symbolic differential equations containing a number of undetermined
coefficient functions, and seeks a solution formula which holds independently
of the special form these coefficients take.

In the first part of [4] he claims to show that a symbolic first-order
ordinary differential equation can be integrated by quadratures in this
sense if and only if it arises from the linear first-order equation $(1)$
by means of a transformation of the unknown variable $y$.
In the second part he proceeds to find criteria for a given equation to
have this property, and concludes, among other things, that the linear
second-order equation (which is intimately connected with the non-linear
first-order Riccati equation) cannot be integrated by quadratures in general.

Despite Maximovi\v c's professed intention `to lay the foundation of an
entirely new theory' of an importance `comparable to that of the theory of
general algebraic equations' ([4], p.I), his work does not seem to have
found a wider audience.
Probably the simple outward reason for this is its practical inaccessibility
as a monograph printed at the Imperial University at Kazan', unavailable even
at the larger American and Western European libraries; thus Ritt in his
classical treatise [7], p.77, regrets that as he `has not been able to
secure Maximovich's paper or any account of it except those given in an
abstract in the {\it Jahrbuch\/} and in one in the Paris {\it Comptes Rendus\/},
he is unable to make a definite statement in regard to it.'
When mentioning Maximovi\v c's result in his annotations to Lie's Collected
Works ([2], p.686), Engel also refers to Vasil'ev's review in the
{\it Jahrbuch\/} [8] only.
Moreover, Maximovi\v c's statements are often obscure and open to
interpretation.
Nevertheless, his work seems to contain some original, useful and justifiable
ideas which deserve to be brought to light.

It is the purpose of the present paper to give a precise definition of
integrability by quadratures along the lines of Maximovi\v c, and to prove
that first-order equations which can be integrated in this sense are
essentially linear.
Although [4] has served as a source of inspiration, no attempt is made to
reconstruct Maximovi\v c's thought.
Specifically, our exposition fundamentally differs from his in that
we consider, instead of a symbolic ordinary differential equation, i.e.
a class of ordinary differential equations of similar structure, a single
differential equation without any hypotheses regarding its structure or
coefficients.
This clearly does not restrict the generality of our main result (Cor.
1.7) as compared to the corresponding theorem of Maximovi\v c.
We assume however (as does Maximovi\v c implicitly) that all initial value
problems have unique solutions defined on a fixed interval.

The paper is organized as follows.
In Section 1 we give the definition of integrability by quadratures, and
state the main result that, up to diffeomorphic transformations of the
unknown variable, the linear first-order equation is the only equation
which can be integrated by quadratures (Cor. 1.7).
This result is a consequence of a normal-form theorem (Thm. 1.6) for
integrals by quadratures which represent the general solution of a first-order
ordinary differential equation. This theorem states that the number of
quadratures can always be reduced to two, which moreover enter the integral
in a very specific way.
Sections 2 and 3 are devoted to the proof of Thm. 1.6.
In Section 2 a property (Lemma 2.1) is derived from the technical requirement
of `independence' of the quadratures, which is fundamental to the actual
reduction procedure described in Section 3.
We remark that although the functions occurring in the quadrature expression
must be assumed to be elementary in some sense in order to avoid a
tautology (cf. Remark 3 to Def. 1.3), the proof of Thm 1.6 does not make
use of the fact or nature of this elementarity, only certain mild
regularity properties are assumed.
In an appendix, we apply Cor. 1.7 to show that the linear second-order
equation can be integrated by quadratures in a somewhat restricted sense
(which excludes, in particular, the well-known examples by D Bernoulli, cf.
[7] VI \S 3) only if it has constant coefficients.
We conclude with a brief account of Maximovi\v c's life, based on the
information given in the preface of [4] and in the obituary notice [9].

\bigskip
{\bf 1 Effectively one-parametric integrals by quadratures.}

\bigskip
In this section we give a definition of integrability by quadratures modelled
roughly on [4], and state (Thm. 1.6) that if such an integral is essentially
one-parametric, e.g. if it represents the general solution of a sufficiently
well-behaved first-order ordinary differential equation, then it can be
reduced to a simple normal form which is very similar to the integral of the
linear equation $(2)$.
From this our main result (Cor. 1.7) follows.

Throughout the paper, we fix an interval $I \subset \R$,
and a point $x_0 \in I$.

\medskip
{\bf 1.1 Definition.}\qquad
We define a {\it system of quadratures\/} recursively as follows:

(i) If $\varphi_1 : I\rightarrow \R$ is a locally integrable function, then,
with
$$
  s_1(x) := \int_{x_0}^x \varphi_1(t)\,dt
  \qquad (x\in I),
$$
$(s_1)$ is a system of quadratures.

(ii) If $(s_1,\dots, s_n)$ is a system of quadratures, $n\in\N$, and
$\varphi_{n+1} : I \times \R^n \rightarrow \R$ is locally integrable in the
first and continuously differentiable in the other variables, then setting
for $x\in I, c_1, \dots, c_n \in \R$,
$$
  s_{n+1}(x, c_1, \dots, c_n) := \int_{x_0}^x \varphi_{n+1}(t, s_1(t) + c_1,
  \dots, s_n(t, c_1, \dots, c_{n-1}) + c_n)\,dt,
$$
$(s_1, \dots, s_n, s_{n+1})$ is a system of quadratures.
We call $\varphi_j$ the {\it integrand\/} of $s_j$.

\medskip
{\it Remark.}\qquad
Thus a system of quadratures is a collection of finite quadrature expressions
which is {\it ordered\/} in that a quadrature can occur, along with its
integration constant, only in the integrands of quadratures with higher index;
and {\it complete\/} in that all quadratures occurring in integrands are
included in the system.
Note that any collection of quadratures which is complete in the latter
sense can be re-arranged to form a system of quadratures, since $\varphi_n$
is permitted to be a constant function of one or more of its arguments.

\medskip
{\bf 1.2 Definition.}\qquad
A system $(s_1, \dots, s_n)$ of quadratures is called {\it independent\/}
if there exist $c_1, \dots, c_{n-1} \in\R$ such that the functions
$\varphi_1, \varphi_2(\Any, c_1), \dots, \varphi_n(\Any, c_1, \dots, c_{n-1})
: I\rightarrow\R$ are linearly independent ($\varphi_j$ being the integrand
of $s_j$, $j\in\{1,\dots,n\}$).

\medskip
{\bf 1.3 Definition.}\qquad
A function $F : \R^n \rightarrow \R$, $n\in\N$, is called {\it admissible\/}
if it is four times continuously differentiable, and $\partial_n F$
has no zeros.

Let ${\it\Theta} : I \times \R \rightarrow \R$ be continuously differentiable
in the second variable such that $\partial_2{\it\Theta}$ has no zeros,
$(s_1, \dots, s_n)$ an independent system of quadratures,
and $F : \R^n \rightarrow \R$ an admissible function, $n\in\N$.
Then we call the family of functions
$$
  f(x, c_1, \dots, c_n) = {\it\Theta}(x, F(s_1(x) + c_1,
  \dots, s_n(x, c_1, \dots, c_{n-1}) + c_n)) \qquad (x\in I, c\in\R^n),
$$
an {\it integral by quadratures}.

\medskip
{\it Remarks.}\qquad
1. The name `integral' reflects that we are primarily interested in families
$f(\Any, c_1, \dots, c_n)$ which represent the general solution of an
ordinary differential equation.

2. Maximovi\v c claims that the above structure of an integral by quadratures
replaces, without loss of generality, the more general expression
$$
  f(x, c_1, \dots, c_n) = F(x, s_1(x) + c_1, s_2(x, c_1) + c_2, \dots,
  s_n(x, c_1, \dots, c_{n-1}) + c_n)
$$
if $f$ represents the general solution of a first-order ordinary differential
equation ([4] \S II).

3. Note that according to the above definition,
{\it every\/} one-parameter family of functions (with
sufficiently regular dependence on the parameter) can be represented by
an integral by quadratures: given $f(\Any, c) : I \rightarrow \R$,
$c \in \R$, we take any locally integrable function $\varphi$ and set
$s(x) := \int_{x_0}^x \varphi$, ${\it\Theta}(x, c) := f(x, c - s(x))$,
$F(c) := c$ $(x\in I, c\in \R)$;
then trivially
$$
  f(x, c) = {\it\Theta}(x, F(s(x) + c))
\qquad\qquad (x\in I, c\in \R).
$$
In order to give some meaning to the concept of `integrability by quadratures',
it is therefore necessary to add the assumption that ${\it\Theta}$ and / or
$F$ be elementary in some sense.
We emphasize, however, that neither the fact nor the nature of this
elementarity are made use of in the considerations of this paper, and we
therefore do not refer to it specifically in our definition.

4. The assumptions that the system of quadratures be independent, 
and the conditions on ${\it\Theta}$ and $F$ are of a
technical nature, necessary in the arguments of Section 3.

\medskip
{\bf 1.4 Definition.} (cf. [3] \S II, D\'ef. I)\qquad
Let $n, m \in \R$.
Two families $f(\Any, c), c\in\R^n$, and
$g(\Any, d), d\in \R^m$, of functions $: I \rightarrow \R$ are called
{\it equivalent\/} if for each $c\in\R^n$ there is a $d\in\R^m$ such that
$f(\Any, c) \equiv g(\Any, d)$, and vice versa.

\medskip
Now consider a first-order ordinary differential equation which, for each
real initial value at $x_0$, has a unique solution defined at least on $I$.
We say that such an equation can be {\it integrated by quadratures\/}
if its general solution on $I \times \R$ is equivalent to an integral by
quadratures.
As the general solution is a one-parameter family of functions parametrized
by their value at $x_0$, the integral, though containing $n$ free integration
constants, has only one effective parameter, which means that the integration
constants are not distinct, but compensate each other.
This observation is central to Maximovi\v c's work.

\medskip
{\bf 1.5 Definition.} \qquad
A family of functions $f(\Any, c) : I \rightarrow \R$, $c\in \R^n$,
is called {\it effectively one-parametric\/}
if it is equivalent to a one-parameter family
$\hat f(\Any, d) : I \rightarrow \R$, $d\in \R$,
with
$\hat f(x_0, d) = d$ $(d\in \R)$.

\medskip
{\it Example.} \qquad
Let $p, q : \R \rightarrow \R$ be linearly independent, locally integrable
functions; then with $\varphi_1(x) := p(x)$,
$\varphi_2(x, c_1) := q(x)\,e^{c_1}$
$(x, c_1 \in \R)$, $(s_1, s_2)$ is an independent system of quadratures.
If ${\it\Theta} : \R^2 \rightarrow \R$ is continuously differentiable in the
second variable, $\partial_2{\it\Theta} \neq 0$,
and ${\it\Theta}(0, \Any)$ is surjective,
then
$$
  f(x, c_1, c_2) := {\it\Theta}\left(x, e^{-(\int_0^x p) - c_1}\,
  \left(\int_0^x q(t)\,e^{(\int_0^t p) + c_1}\,dt + c_2\right)\right)
\qquad
  (x, c_1, c_2 \in \R)
$$
is an effectively one-parametric integral by quadratures
(with $F(c_1, c_2) = e^{-c_1} c_2$).

Indeed, it is easy to see that the parameter $c_1$ can be eliminated, and
$f(\Any, c_1, c_2)$ is equivalent to
$\hat f(\Any, d) := f(\Any, 0, \Theta^{-1}(0, d))$
(where $\Theta^{-1}(0, \Any)$ is the inverse of $\Theta(0, \Any)$).

\medskip
Our main result shows that the structure of $\hat f$ in this simple
example is in fact universal; specifically,
we are going to prove

\medskip
{\bf 1.6 Theorem.} (Normal form of effectively one-parametric integrals by
quadratures.)\qquad
{\it
If $f$ is an effectively one-parametric integral by quadratures,
then there are locally integrable functions $p, q : I \rightarrow \R$
and a function ${\it\hat\Theta} : I \times\R\rightarrow\R$
which is continuously differentiable in the second variable, such that
$f$ is equivalent to the family
$$
  \hat f(x, C) := {\it\hat\Theta}\left(x, \exp\left(-\int_{x_0}^x p\right)\,
   \left(\int_{x_0}^x q(t)\,\exp\left(\int_{x_0}^t p\right)\,dt + C
   \right)\right),
$$
$(x\in I, C\in\R)$.
}

\medskip
{\it Remark.}\qquad
Moreover, if $f$ is given in the form
$$
  f(x, c_1, \dots, c_n)
  = {\it\Theta}(x, F(s_1(x) + c_1, \dots, s_n(x, c_1, \dots, c_{n-1}) + c_n)),
$$
then
--- as will be apparent from the proof of Thm. 1.6 ---
${\it\hat\Theta}(x, c) = {\it\Theta}(x, F(0, \dots, 0, c))$
$(x\in I, c\in\R)$.
The final $p$ and $q$ arise from $F$ and the integrands of $s_1, \dots, s_n$
by the basic arithmetical operations, exponentiation, differentiation and
integration.
Thus the process leading to the normal form is completely transparent,
and in particular the elementarity of ${\it\hat\Theta}$ could be inferred
immediately from that of ${\it\Theta}$ and $F$.

\medskip
{\bf 1.7 Corollary.}\qquad
{\it
A first-order ordinary differential equation which,
for each real initial value at $x_0$,
has a unique solution which is defined at least on $I$,
can be integrated by quadratures only if
it arises from the linear first-order equation
by a diffeomorphic transformation of the unknown variable.
}

\medskip
{\it Remarks.}\qquad
1. As a consequence,
a first-order differential equation
which can be integrated by quadratures
is of the type
$$
  \partial_1 \Phi(x, y) + \partial_2 \Phi(x, y) \cdot y'
  = q(x) - p(x)\,\Phi(x, y),
$$
with $p, q : I \rightarrow \R$ locally integrable
and $\Phi : I \times \R \rightarrow \R$ continuously differentiable,
$\partial_2 \Phi \neq 0$.
($\Phi$ is the inverse of $\hat \Theta$ with respect to the second variable.)

If $p \equiv 0$, $\partial_1 \Phi \equiv 0$, then the resulting equation
$\Phi'(y)\,y' = q(x)$ has separated variables;
note that its integral by quadratures,
$$
  f(x, c) = \Phi^{-1}\left(\int_{x_0}^x q(t)\,dt + c \right),
$$
contains only one quadrature.
The second quadrature in the usual solution formula
for differential equations of separated type,
which arises in the process of passing from $\Phi'$ to $\Phi^{-1}$,
does not introduce an independent integration constant
and does not figure explicitly in Maximovi\v c's form of the integral.

2. By the Picard-Lindel\"of Theorem,
differential equations of the type
$$
  y' = g(y, x)
$$
satisfy the requirement of having a unique solution on $I$ for all initial data
if $g$ is uniformly Lipschitz in the first variable.
Our hypothesis excludes cases with movable singularities,
such as the Riccati equation.
However,
the Riccati equation $y' = y^2+ Q$ can easily be transformed into
the Pr\"ufer equation $\vartheta' = \cos^2 \vartheta + Q \sin^2 \vartheta$,
which is uniformly Lipschitz in $\vartheta$.

It seems likely that the statement and proof of Thm. 1.6 can be made local and
that this restriction can be removed in this way;
yet we prefer to study the global general solution in order to avoid
notational complications which would obscure the underlying ideas.

\bigskip
{\bf 2 Independent systems of quadratures.}

\bigskip
In this section
we prove the following property of independent systems of quadratures
which will play a central role in the reduction process of Section 3.

\medskip
{\bf 2.1 Lemma.}\qquad
{\it
Let $(s_1, \dots, s_n)$ be an independent system of quadratures.
If $g : \R^n\rightarrow \R$ is a continuously differentiable function such
that
$$
  g(s_1(x) + c_1, \dots, s_n(x, c_1, \dots, c_{n-1}) + c_n)
  = g(c_1, \dots, c_n)
\eqno (3)
$$
$(x\in I, c\in\R^n)$, then $g$ is constant.
}

\medskip
{\it Remark.}\qquad
This property could also be used as a definition of independence for the
purposes of Section 3. However, our definition emphasizes the local character
of independence, and proves more convenient in the reduction process.

It may seem surprising that the global property stated in the above assertion
follows from the very local condition of linear independence of the integrands
at {\it one} point; as will be apparent from the proof, the
decisive point is the hierarchical order of the quadratures which ensures
that the $j$-th component of the curve
$$
  t \mapsto (s_1(t) + c_1, \dots, s_n(t, c_1, \dots, c_{n-1}) + c_n)
  \qquad (t\in I)
$$
--- along which $g$ is constant by hypothesis --- depends only on $c_1, \dots,
c_{j-1}$, but not on $c_j$.

\medskip
{\it Proof.}\qquad
The independence of the quadratures means that there is a point
$(\hat c_1, \dots, \hat c_{n-1}) \in \R^{n-1}$ such that the $n$ functions
$$
  \varphi_1, \varphi_2(\Any, \hat c_1), \dots, \varphi_n(\Any, \hat c_1, \dots,
  \hat c_{n-1})
$$
are linearly independent.
It is easy to prove by induction over $k$ that then there exist $n$ points
$x_1, \dots, x_n \in I$ such that the matrix
$ \left(\varphi_j(x_i, \hat c_1, \dots, \hat c_{j-1})\right)_{i,j\in
   \{1,\dots,k\}} $
has rank $k$, for $k\in\{1, \dots, n\}$.

Let $\gamma_1, \dots, \gamma_n \in \R$ and $\xi\in I$.
In the equation $(3)$ set
$$\eqalign{
  c_1 &:= \gamma_1 - s_1(\xi),
\cr
  c_2 &:= \gamma_2 - s_2(\xi, \gamma_1 - s_1(\xi)),
\cr
  \dots, c_n &:= \gamma_n - s_n(\xi, \gamma_1 - s_1(\xi), \dots,
        \gamma_{n-1} - s_{n-1}(\xi, \dots)).
\cr }$$

Then, differentiating with respect to $x$ at $x = \xi$, we obtain
$$
  0 = \sum_{j=1}^n \partial_j g(\gamma_1, \dots, \gamma_n)\,
      \varphi_j(\xi, \gamma_1, \dots, \gamma_{j-1});
\eqno (4)
$$
this holds for all $\xi\in I$ and $\gamma_1, \dots, \gamma_n \in \R$.

Now fix $\hat c_n \in \R$ arbitrary, and set
$\hat g := g(\hat c_1, \dots, \hat c_n)$
(with $\hat c_1, \dots, \hat c_{n-1}$ as above).
We shall prove by induction over $k\in \{0, \dots, n\}$ that
$$
  g(\hat c_1, \dots, \hat c_{n-k}, c_{n-k+1}, \dots, c_n) = \hat g
  \qquad (c_{n-k+1}, \dots, c_n \in \R).
$$
In the case $k=0$ this is clear from the definition of $\hat g$.

Now assume we know the assertion is true for some $k\in\{0, \dots, n-1\}$;
we want to show that it also holds for $k+1$.

As the vectors $(\varphi_1(x_i), \varphi_2(x_i, \hat c_1), \dots,
  \varphi_{n-k}(x_i, \hat c_1, \dots, \hat c_{n-k-1}))_{i\in\{1, \dots, n-k\}}$
are linearly independent, there are constants $\alpha_1, \dots, \alpha_{n-k}
\in \R$ such that
$$
  \sum_{i=1}^{n-k} \alpha_i\, \varphi_j(x_i, \hat c_1, \dots, \hat c_{j-1})
  = \delta_{j, n-k}
  \qquad (j\in \{1, \dots, n-k\}).
$$
Multiplying $(4)$ by $\alpha_i$ and summing over $i\in\{1, \dots, n-k\}$, and
setting for

$j\in\{n-k+1, \dots, n\}$
$$
  \psi_j(c_{n-k}, \dots, c_{j-1}) := 
       \sum_{i=1}^{n-k} \alpha_i\, \varphi_j(x_i, \hat c_1, \dots,
          \hat c_{n-k-1}, c_{n-k}, \dots, c_{j-1}),
$$
we obtain the linear first-order partial differential equation
$$\eqalign{
  0 &= \partial_{n-k} g(\hat c_1, \dots, \hat c_{n-k-1}, c_{n-k}, \dots, c_n)
\cr
    &\qquad + \sum_{j=n-k+1}^n \partial_j g(\hat c_1, \dots, \hat c_{n-k-1},
         c_{n-k}, \dots, c_n)\, \psi_j(c_{n-k}, \dots, c_{j-1}),
\cr }$$
with initial data $\hat g$ on the (non-characteristic) surface
$\{(\hat c_1, \dots, \hat c_{n-k})\}\times \R^k$.
In order to solve this equation we consider the characteristic initial value
problems
(cf. [1] I \S 5)
$$\eqalign{
  &\gamma_{n-k}'(t) = 1, \qquad \gamma_{n-k}(0) = \hat c_{n-k},
\cr
  &\gamma_j'(t) = \psi_j(\gamma_{n-k}(t), \dots, \gamma_{j-1}(t)), \qquad
    \gamma_j(0) = \hat\gamma_j \qquad (j\in\{n-k+1, \dots, n\}),
\cr } $$
with $\hat\gamma_{n-k+1}, \dots, \hat\gamma_n \in \R$ arbitrary.

This system of $k+1$ ordinary differential equations fully decouples; the
equations can be solved one after the other by simple integration:
$$\eqalign{
  \gamma_{n-k}(t) &= \hat c_{n-k} + t,
\cr
  \gamma_{n-k+1}(t) &= \hat\gamma_{n-k+1}
        + \int_0^t \psi_{n-k+1}(\gamma_{n-k}(s))\,ds,
\cr
  \dots, \gamma_n(t) &= \hat\gamma_n
        + \int_0^t \psi_n(\gamma_{n-k}(s), \dots, \gamma_{n-1}(s))\,ds.
\cr }$$

$g$ is constant along the curves
$ t \mapsto (\hat c_1, \dots, \hat c_{n-k-1}, \gamma_{n-k}(t), \dots,
   \gamma_n(t)), \quad t\in\R;$
thus
$$
 g(\hat c_1, \dots, \hat c_{n-k-1}, \hat c_{n-k} + t, \gamma_{n-k+1}(t),
    \dots, \gamma_n(t))
 = g(\hat c_1, \dots, \hat c_{n-k-1}, \hat c_{n-k}, \hat\gamma_{n-k+1}, \dots,
    \hat\gamma_n).
$$
Now given $c_{n-k}, \dots, c_n \in \R$, we set
$$\eqalign{
  t &:= c_{n-k} - \hat c_{n-k},
\cr
  \hat\gamma_{n-k+1} &:= c_{n-k+1} - \int_0^t \psi_{n-k+1}(\gamma_{n-k}(s))\,ds,
\cr
  \dots,\hat\gamma_n &:= c_n - \int_0^t \psi_n(\gamma_{n-k}(s), \dots,
    \gamma_{n-1}(s))\,ds,
\cr} $$
and find
$$
 g(\hat c_1, \dots, \hat c_{n-k-1}, c_{n-k}, c_{n-k+1}, \dots, c_n)
 = g(\hat c_1, \dots, \hat c_{n-k-1}, \hat c_{n-k},
     \hat\gamma_{n-k+1}, \dots, \hat\gamma_n)
 = \hat g.  
\eqno \square
$$

\bigskip
{\bf 3 The reduction procedure.}

\bigskip
In this section we prove Thm. 1.6.
As a starting point
we show that the fact that
the integral by quadratures is essentially one-parametric
can be expressed by a simple relationship
between its partial derivatives with respect to the free integration constants
(Prop. 3.1).
Then we provide two auxiliary propositions
which capture two fundamental processes
which are repeatedly applied in the reduction procedure:
the explicit solution of a very simply-structured first-order linear partial
 differential equation (Lemma 3.2);
and the inclusion of an arbitrary function of a number of quadratures
in the integrand of a quadrature with higher index (Lemma 3.3).
Then we proceed to the heart of the matter in Props. 3.4 and 3.5,
which demonstrate how the number of quadratures in the integral
can be reduced iteratively
by combining two quadratures into one with the help of Prop. 3.1.
After the proof of Prop. 3.5
we outline the reduction algorithm,
and finish the proof of Thm. 1.6.

\medskip
{\bf 3.1 Proposition.}\qquad
{\it
If $f(\Any, c) : I \rightarrow \R$, $c \in \R^n$,
is an effectively one-parametric family
continuously differentiable with respect to $c$,
then the following {\rm Fundamental Equality} ({\it for the pair $c_i, c_j$}),
$i,j\in\{1,\dots,n\}$, holds:
$$
  \partial_{1+i} f(x, c)\, \partial_{1+j} f(x_0, c)
= \partial_{1+j} f(x, c)\, \partial_{1+i} f(x_0, c)
\qquad (x\in I, c\in\R^n).
$$
}

\smallskip
{\it Proof.}\qquad
By the equivalence of $f$ and
a one-parameter family $\hat f$ parametrized by its values at $x_0$,
there is a function $d : \R^n \rightarrow \R$ such that
$f(\Any, c) = \hat f(\Any, d(c))$ \quad $(c\in\R^n)$;
in particular $f(x_0, c) = \hat f(x_0, d(c)) = d(c)$ \quad $(c\in\R^n)$.
Differentiating the identity
$f(x, c_1, \dots, c_n) = \hat f(x, f(x_0, c_1, \dots, c_n))$
with respect to $c_j$, we find
$$
  \partial_{1+j} f(x, c_1, \dots, c_n)
  = \partial_2 \hat f(x, f(x_0, c_1, \dots, c_n))\,
    \partial_{1+j} f(x_0, c_1, \dots, c_n)
$$
$(x\in I, c_1, \dots, c_n \in \R)$,
from which the Fundamental Equality follows
upon elimination of $\partial_2 \hat f$.
\hfill $\square$

\medskip
{\bf 3.2 Lemma.}\qquad
{\it
Let $H : \R^n \rightarrow \R$ be continuously differentiable, and
$a, b : \R^{n-1} \rightarrow \R$ continuous such that
$$
  0 = \partial_{n-1} H(x_1, \dots, x_n) - (a(x_1, \dots, x_{n-1})\, x_n
      + b(x_1, \dots, x_{n-1}))\, \partial_n H(x_1, \dots, x_n).
$$
Then, with $G(x_1, \dots, x_{n-2}, x_n) := H(x_1, \dots, x_{n-2}, 0, x_n)$,
$$\eqalign{
  H(x_1, \dots, x_n) = G\Big(x_1,& \dots, x_{n-2},
   \exp\Big(\int_0^{x_{n-1}} a(x_1, \dots, x_{n-2}, t)\,dt\Big)\,x_n
\cr
   &+ \int_0^{x_{n-1}} b(x_1, \dots, x_{n-2}, t)\,
   \exp\Big(\int_0^t a(x_1, \dots, x_{n-2}, s)\,ds\Big)\,dt\Big)
\cr }$$
$(x_1, \dots, x_n \in \R)$.
}

\smallskip
{\it Proof.}\qquad
We apply the standard procedure to solve a (quasi-) linear partial differential
equation of first order, cf. [1] I \S 5.
For fixed $x_1, \dots, x_{n-2} \in \R$,
we solve the characteristic initial value problems
$$\eqalign{
  \xi_{n-1}'(t) &= 1, \qquad \xi_{n-1}(0) = 0,
\cr
  \xi_n'(t) &= -a(x_1, \dots, x_{n-2}, \xi_{n-1}(t))\,\xi_n(t)
               - b(x_1, \dots, x_{n-2}, \xi_{n-1}(t)),
  \qquad \xi_n(0) = \hat\xi_n\in\R.
\cr}$$
The solution is given by formula $(2)$:
$$\eqalign{
  \xi_{n-1}(t) &= t,
\cr
  \xi_n(t) &= \exp\left(-\int_0^t a(x_1, \dots, x_{n-2}, s)\, ds \right)\,
\cr
     &\qquad\times
           \left(\hat\xi_n - \int_0^t b(x_1, \dots, x_{n-2}, s)\,
            \exp\left(\int_0^t a(x_1, \dots, x_{n-2}, r)\, dr \right)\,ds\right)
\cr}$$
$(t \in \R)$.
Now setting for $x_{n-1}, x_n \in \R$
$$\eqalign{
  \hat\xi_n
   := &\exp\left(\int_0^{x_{n-1}} a(x_1, \dots, x_{n-2}, t)\, dt\right)\, x_n
\cr
   &+ \int_0^{x_{n-1}} b(x_1, \dots, x_{n-2}, t)\,
      \exp\left(\int_0^t a(x_1, \dots, x_{n-2}, s)\, ds\right)\, dt
\cr}$$
and observing that
$$
  0 = {d \over dt} H (x_1, \dots, x_{n-2}, \xi_{n-1}(t), \xi_n(t))
  \qquad (t \in \R),
$$
the assertion follows.
\hfill $\square$

\medskip
{\bf 3.3 Lemma.}\qquad
{\it
Let $(s_1, \dots, s_n)$ be an independent system of quadratures, $n \ge 2$,
and $B : \R^{n-1} \rightarrow \R$  twice continuously differentiable.
Then there is a quadrature $\hat s_n$ such that $(s_1, \dots, s_{n-1},
\hat s_n)$ is an independent system of quadratures, and
$$\eqalign{
  s_n(x, c_1, \dots, c_{n-1}) + B&(s_1(x) + c_1, \dots, s_{n-1}(x, c_1, \dots,
   c_{n-2}) + c_{n-1})
\cr
  &= \hat s_n(x, c_1, \dots, c_{n-1}) + B(c_1, \dots, c_{n-1})
\cr}$$
$(x\in I, c_1, \dots, c_{n-1} \in \R)$.
}

\smallskip
{\it Proof.}\qquad
Let $\varphi_1, \dots, \varphi_n$ denote the integrands of $s_1, \dots, s_n$.
Set
$$\eqalign{
  \hat\varphi_n(x, c_1,& \dots, c_{n-1})
  := \varphi_n(x, c_1, \dots, c_{n-1})
   + \sum_{j=1}^{n-1} \partial_j B(c_1, \dots, c_{n-1})\,
       \varphi_j(x, c_1, \dots, c_{j-1});
\cr
  \hat s_n(x, c_1,& \dots, c_{n-1}) := \int_{x_0}^x \hat\varphi_n(t,
    s_1(t) + c_1, \dots, s_{n-1}(t, c_1, \dots, c_{n-2}) + c_{n-1})\,dt
\cr}$$
$(x\in I, c_1, \dots, c_n \in \R)$.
\hfill $\square$

\medskip
In the following the last (highest-order) quadrature plays a special role and
is therefore distinguished in the notation.

\medskip
{\bf 3.4 Proposition.}\qquad
{\it
Let $f$ be an effectively one-parametric integral by quadratures of the form
$$
  f(x, c_1, \dots, c_n, C) := {\it\Theta}(x,
    F(s_1(x) + c_1, \dots, s_n(x, c_1, \dots,
    c_{n-1}) + c_n, S(x, c_1, \dots, c_n) + C))
$$
$(x\in I, c_1, \dots, c_n, C \in \R)$, $n\in\N$.

Set $G(c_1, \dots, c_{n-1}, C) := F(c_1, \dots, c_{n-1}, 0, C)$\qquad
$(c_1, \dots, c_{n-1}, C \in \R)$.

Then one of the following statements is true:

1. There is a quadrature $\hat S$ such that
$(s_1, \dots, s_{n-1}, \hat S)$ is an independent system of quadratures,
and the integral by quadratures
$$\eqalign{
  \hat f(x, &c_1, \dots, c_{n-1}, \hat C)
\cr
  &:= {\it\Theta}(x,
      G(s_1(x) + c_1, \dots, s_{n-1}(x, c_1, \dots, c_{n-2}) + c_{n-1},
             \hat S(x, c_1, \dots, c_{n-1}) + \hat C))
\cr}$$
$(x\in I, c_1, \dots, c_{n-1}, \hat C \in \R)$,
is equivalent to $f$.

2. There are quadratures $s, \hat S$ such that
$(s_1, \dots, s_{n-1}, s, \hat S)$ is an independent system of quadratures,
and the integral by quadratures
$$\eqalign{
  \hat f(x, c_1, \dots, c_{n-1}, c, \hat C)
  := {\it\Theta}(x,
       G&(s_1(x) + c_1, \dots, s_{n-1}(x, c_1, \dots, c_{n-2}) + c_{n-1},
\cr
  & e^{s(x, c_1, \dots, c_{n-1}) + c}\,(\hat S(x, c_1, \dots, c_{n-1}, c)
     + \hat C)))
\cr}$$
$(x\in I, c_1, \dots, c_{n-1}, c, \hat C \in \R)$,
is equivalent to $f$.
}

\medskip
{\it Proof.}\qquad
Writing out the Fundamental Equality for $f$ and the pair $c_n, C$
in terms of the above expression,
and dividing by factors of $\partial_2{\it\Theta}$, $\partial_{n+1} F \neq 0$,
we find, for all $x\in I, c_1, \dots, c_n, C \in \R$:
$$
  {\partial_n F \over \partial_{n+1} F} (c_1, \dots, c_n, C)
  = {\partial_n F \over \partial_{n+1} F} (s_1(x) + c_1, \dots,
      S(x, c_1, \dots, c_n) + C)
   + \partial_{1+n} S(x, c_1, \dots, c_n).
$$
Differentiation of this identity with respect to $C$ yields
$$
  \partial_{n+1}\,{\partial_n F \over \partial_{n+1} F} (c_1, \dots, c_n, C)
  = \partial_{n+1}\,{\partial_n F \over \partial_{n+1} F} (s_1(x) + c_1, \dots,
      S(x, c_1, \dots, c_n) + C).
$$
By the independence of the quadratures and Lemma 2.1 it follows that there is
a constant $\alpha\in\R$ such that
$$
  \partial_{n+1}\,{\partial_n F \over \partial_{n+1} F} (c_1, \dots, c_n, C)
  = \alpha
\qquad (c_1, \dots, c_n, C \in \R);
$$
and consequently
$$
  \partial_n F(c_1, \dots, c_n, C)
  = (\alpha\, C + \beta(c_1, \dots, c_n))\,\partial_{n+1} F(c_1, \dots, c_n, C)
\qquad (c_1, \dots, c_n, C \in \R)
$$
with some twice continuously differentiable function
$\beta : \R^n \rightarrow \R$.
Hence by Lemma 3.2
$$
  F(c_1, \dots, c_n, C)
  = G\Big(c_1, \dots, c_{n-1},
      e^{\alpha c_n} \Big(C + e^{-\alpha c_n}
      \int_0^{c_n} \beta(c_1, \dots, c_{n-1}, t)\, e^{\alpha t}\, dt\Big)\Big)
$$
$(c_1, \dots, c_n, C \in \R)$.

Now there are two different cases to consider:

{\bf 1st case:} $\alpha = 0$.
\qquad
Then
$$
  F(c_1, \dots, c_n, C) = G(c_1, \dots, c_{n-1}, C + B(c_1, \dots, c_n)),
$$
$$
  B(c_1, \dots, c_n) := \int_0^{c_n} \beta(c_1, \dots, c_{n-1}, t)\, dt
\qquad (c_1, \dots, c_n \in \R).
\leqno \hbox{where}
$$
By Lemma 3.3 there is a quadrature $S'$ (here and in corresponding cases below
the prime is used merely for distinction and does not signify a derivative)
such that
$(s_1, \dots, s_n, S')$ is an independent system of quadratures and
$$\eqalign{
  f(x, c_1, \dots, c_n, C)
  = {\it\Theta}(x,
     G(s_1(x) + c_1, &\dots, s_{n-1}(x, c_1, \dots, c_{n-2}) + c_{n-1},
\cr
  & S'(x, c_1, \dots, c_n) + C + B(c_1, \dots, c_n))).
\cr}$$
$(x\in I, c_1, \dots, c_n, C \in \R)$.

Note that the quadrature $s_n$ does no longer occur explicitly as an
argument of $G$,
but only in the integrand of $S'$;
we now show that it is redundant and can be eliminated.

Clearly $\hat C := C + B(c_1, \dots, c_n)$
takes the role of a new integration constant for $S'$.
Thus $f$ is equivalent to the integral by quadratures
$$\eqalign{
  f'(x, c_1, &\dots, c_n, \hat C)
\cr
  &= {\it\Theta}(x,
      G(s_1(x) + c_1, \dots, s_{n-1}(x, c_1, \dots, c_{n-2}) + c_{n-1},
      S'(x, c_1, \dots, c_n) + \hat C))
\cr}$$
$(x\in I, c_1, \dots, c_n, \hat C \in \R)$.

Therefore --- since equivalence is transitive ---
$f'$ is effectively one-parametric, and by Prop. 3.1
the Fundamental Equality holds for the pair $c_n, \hat C$;
its left-hand part vanishes
as the initial value $f'(x_0, c_1, \dots, c_n, \hat C)$ is independent
  of $c_n$:
$$\eqalign{
  0
 = \partial_n G(s_1(x) &+ c_1, \dots, s_{n-1}(x, c_1, \dots, c_{n-2}) + c_{n-1},    S'(x, c_1, \dots, c_n) + \hat C)
\cr
  &\times \partial_{1+n} S'(x, c_1, \dots, c_n)
   \, \partial_n G(c_1, \dots, c_{n-1}, \hat C);
\cr}$$
since $\partial_n G \neq 0$, it follows that
$\partial_{1+n} S' = 0$.
Differentiating with respect to $x$
(and denoting the integrand of $S'$ by ${\it\Phi}'$)
we find
$$
  0
  = \partial_{1+n} {\it\Phi}'(x, s_1(x) + c_1, \dots, s_n(x, c_1, \dots,
      c_{n-1}) + c_n)
\qquad (x\in I, c_1, \dots, c_n \in \R).
$$
This means that for fixed $x, c_1, \dots, c_{n-1}$ the integrand has the
same value for all real values of $c_n$; therefore setting
${\it\hat\Phi}(x, c_1, \dots, c_{n-1})
  := {\it\Phi}'(x, c_1, \dots, c_{n-1}, 0),$
$$
  \hat S(x, c_1, \dots, c_{n-1})
  := \int_{x_0}^x {\it\hat\Phi}(t, s_1(t) + c_1, \dots, s_{n-1}(t, c_1, \dots,
              c_{n-2}) + c_{n-1})\, dt
$$
has the required properties.

{\bf 2nd case:} $\alpha \neq 0$.
\qquad
Denoting by $\varphi_n$ and ${\it\Phi}$ the integrands of $s_n$ and $S$, resp.,
define new quadratures $s$, $S'$
by giving their integrands
$$\eqalign{
  \varphi&(x, c_1, \dots, c_{n-1})
  := \alpha\, \varphi_n(x, c_1, \dots, c_{n-1})
\cr
  {\it\Phi}'&(x, c_1, \dots, c_{n-1}, c)
  := {\it\Phi}(x, c_1, \dots, c_{n-1}, c/\alpha):
\cr
  s&(x, c_1,\dots, c_{n-1})
  := \int_{x_0}^x \varphi(t, s_1(t) + c_1, \dots,
     s_{n-1}(t, c_1, \dots, c_{n-2}) + c_{n-1}))\,dt,
\cr
  S'&(x, c_1,\dots, c_{n-1}, c)
  := \int_{x_0}^x {\it\Phi}'(t, s_1(t) + c_1, \dots,
     s_{n-1}(t, c_1, \dots, c_{n-2}) + c_{n-1},
\cr
    &\qquad\qquad\qquad\qquad \qquad\qquad s(t, c_1, \dots, c_{n-1}) + c)
   \,dt.
\cr}$$

Then $(s_1, \dots, s_{n-1}, s, S')$
is an independent system of quadratures,
and $S(x, c_1, \dots, c_n) = S'(x, c_1, \dots, c_{n-1}, \alpha\,c_n)$.
Furthermore, set
$$
  B(c_1, \dots, c_{n-1}, c)
  := e^{-c} \int_{x_0}^{c/\alpha} \beta(c_1, \dots, c_{n-1}, t)\, e^{\alpha t}
  \,dt
  \qquad (c_1, \dots, c_n, c \in \R);
$$
then by Lemma 3.3 we find a quadrature $\hat S$ such that
$(s_1, \dots, s_{n-1}, s, \hat S)$ is an independent system of quadratures,
and
$$\eqalign{
  f(x, c_1, \dots, c_n, &C)
  = {\it\Theta}(x,
      G(s_1(x) + c_1, \dots, s_{n-1}(x, c_1, \dots, c_{n-2}) + c_{n-1},
\cr
  & e^{s(x, c_1, \dots, c_{n-1}) + \alpha c_n}
     (\hat S(x, c_1, \dots, c_{n-1}, \alpha c_n) + C
      + B(c_1, \dots, c_{n-1}, \alpha c_n)) )).
\cr}$$
Hence,
with $c := \alpha\, c_n$ and $\hat C = C + B(c_1, \dots, c_{n-1}, \alpha c_n)$,
$f$ is equivalent to the integral $\hat f$.
\hfill $\square$

\medskip
{\bf 3.5 Proposition.}\qquad
{\it
Let $f$ be an effectively one-parametric integral by quadratures of the form
$$\eqalign{
  f(x, c_1, \dots, c_n, c, C)
  := {\it\Theta}(x,
       F(s_1(x) + c_1, \dots, & s_n(x, c_1, \dots, c_{n-1}) + c_n,
\cr
    & e^{s(x, c_1, \dots, c_n) + c} (S(x, c_1, \dots, c_n, c) + C)))
\cr}$$
$(x\in I, c_1, \dots, c_n, c, C \in \R)$,
$n\in \N$.

Set $G(c_1, \dots, c_{n-1}, D) := F(c_1, \dots, c_{n-1}, 0, D)$\qquad
$(c_1, \dots, c_{n-1}, D \in \R)$.

Then there are quadratures $\hat s$, $\hat S$
such that
$(s_1, \dots, s_{n-1}, \hat s, \hat S)$ is an independent system of quadratures
and the integral by quadratures
$$\eqalign{
  \hat f(x, c_1, \dots, c_{n-1}, \hat c, \hat C)
  := {\it\Theta}(x,
       G(s_1(x) + c_1, \dots, & s_{n-1}(x, c_1, \dots, c_{n-2}) + c_{n-1},
\cr
    & e^{\hat s(x, c_1, \dots, c_{n-1}) + \hat c}
       (\hat S(x, c_1, \dots, c_{n-1}, \hat c) + \hat C)))
\cr}$$
$(x\in I, c_1, \dots, c_{n-1}, \hat c, \hat C \in \R)$
is equivalent to $f$.
}

{\it Proof.}\qquad
Writing out the Fundamental Equality for $f$ and the pair $c_n, C$
in terms of ${\it\Theta}$ and $F$,
and dividing by factors of $\partial_2{\it\Theta}$, $\partial_{n+1} F \neq 0$,
we find
$$\eqalign{
  &{\partial_n F \over \partial_{n+1} F} (c_1, \dots, c_n, e^c C)
  \, e^{-c}
\cr
 &= {\partial_n F \over \partial_{n+1} F}
     (s_1(x) + c_1, \dots, s_n(x, \dots) + c_n,
            e^{s(x, \dots) + c} (S(x, \dots) + C))
  \, e^{-s(x, c_1, \dots, c_n) - c}
\cr
 &\qquad + (S(x, c_1, \dots, c_n, c) + C)
   \, \partial_{1+n} s(x, c_1, \dots, c_n)
  + \partial_{1+n} S(x, c_1, \dots, c_n, c).
\cr}$$

Differentiating this equality twice with respect to $C$ we obtain
$$\eqalign{
  &\partial_{n+1} \partial_{n+1} {\partial_n F \over \partial_{n+1} F}
    (c_1, \dots, c_n, e^c C)
  \, e^c
\cr
 &= \partial_{n+1} \partial_{n+1} {\partial_n F \over \partial_{n+1} F}
     (s_1(x) + c_1, \dots, s_n(x, \dots) + c_n,
            e^{s(x, \dots) + c} (S(x, \dots) + C))
  \, e^{s(x, c_1, \dots, c_n) + c}.
\cr}$$

Since the quadratures are independent,
this implies by Lemma 2.1 that
there is a constant $\gamma \in \R$ such that
$$
  \partial_{n+1} \partial_{n+1} {\partial_n F \over \partial_{n+1} F}
    (c_1, \dots, c_n, D)
  \, e^c
  = \gamma
  \qquad\qquad (c_1, \dots, c_n, c, D \in \R).
$$
As $c$ is arbitrary and
the first factor on the left hand side does not depend on $c$,
$\gamma$ must be $0$.
Thus there are twice continuously differentiable functions
$\alpha, \beta : \R^n \rightarrow \R$ such that
$$
  \partial_n F (c_1, \dots, c_n, D)
  = (\alpha(c_1, \dots, c_n)\, D + \beta(c_1, \dots, c_n))
    \, \partial_{n+1} F(c_1, \dots, c_n, D)
$$
$(c_1, \dots, c_n, D \in \R)$.

Hence by Lemma 3.2
$$
  F(c_1, \dots, c_n, D)
  = G(c_1, \dots, c_{n-1},
           e^{A(c_1, \dots, c_n)} D + B'(c_1, \dots, c_n)),
\qquad (c_1, \dots, c_n, D \in \R)
$$
where
$$\eqalign{
  A(c_1, \dots, c_n)
 &:= \int_0^{c_n} \alpha(c_1, \dots, c_{n-1}, t)\, dt,
  \qquad \hbox{and}
\cr
  B'(c_1, \dots, c_n)
 &:= \int_0^{c_n} \beta(c_1, \dots, c_{n-1}, t)\,
      e^{A(c_1, \dots, c_{n-1}, t)}\, dt
\cr}$$
$(c_1, \dots, c_n \in \R)$.

By Lemma 3.3
there is a quadrature $s'$ such that
$$
  s(x, c_1, \dots, c_n)
  + A(s_1(x) + c_1, \dots, s_n(x, c_1, \dots, c_{n-1}) + c_n)
  = s'(x, c_1, \dots, c_n)
  + A(c_1, \dots, c_n).
$$
Denoting by ${\it\Phi}$ the integrand of $S$ and
setting
$$
  {\it\Phi}'(x, c_1, \dots, c_n, c')
  := {\it\Phi}(x, c_1, \dots, c_n, c' - A(c_1, \dots, c_n))
$$
we have the new quadrature
$$
  S'(x, c_1, \dots, c_n, c')
  := \int_{x_0}^x {\it\Phi}'(t, s_1(t) + c_1, \dots,
             s_n(t, c_1, \dots, c_{n-1}) + c_n, s'(t, c_1, \dots, c_n) + c')
  \, dt
$$
$(x\in I, c_1, \dots, c_n, c' \in \R)$;
then $(s_1, \dots, s_n, s', S')$ is an independent system of quadratures,
and
$$\eqalign{
  f(x, c_1, \dots, c_n, c, C)
  &= {\it\Theta}(x,
    G(s_1(x) + c_1, \dots, s_{n-1}(x, c_1, \dots, c_{n-2}) + c_{n-1},
\cr
  & e^{s'(x, c_1, \dots, c_n) + c + A(c_1, \dots, c_n)}
   \, (S'(x, c_1, \dots, c_n, c + A(c_1, \dots, c_n)) + C)
\cr
  &\qquad\qquad
    + B'(s_1(x) + c_1, \dots, s_n(x, c_1, \dots, c_{n-1}) + c_n))).
\cr}$$
With $B(c_1, \dots, c_n, c') := e^{-c'} B'(c_1, \dots, c_n)$,
there is by Lemma 3.3 a quadrature $S''$ such that
$(s_1, \dots, s_{n-1}, s', S'')$ is an independent system of quadratures, and
$$\eqalign{
  f(x, c_1, \dots, c_n, c, C)
  &= {\it\Theta}(x,
      G(s_1(x) + c_1, \dots, s_{n-1}(x, c_1, \dots, c_{n-2}) + c_{n-1},
\cr
  & e^{s'(x, c_1, \dots, c_n) + c + A(c_1, \dots, c_n)}
   \, (S''(x, c_1, \dots, c_n, c + A(c_1, \dots, c_n)) + C
\cr
  &\qquad\qquad
    + B(c_1, \dots, c_n, c + A(c_1, \dots, c_n))))),
\cr}$$
which is clearly equivalent to
$$\eqalign{
  f'(x, c_1, \dots, c_n, \hat c, \hat C)
  &= {\it\Theta}(x,
       G(s_1(x) + c_1, \dots, s_{n-1}(x, c_1, \dots, c_{n-2}) + c_{n-1},
\cr
  &\qquad e^{s'(x, c_1, \dots, c_n) + \hat c}
   \, (S''(x, c_1, \dots, c_n, \hat c) + \hat C))).
\cr}$$

Note that the quadrature $s_n$ occurs only in the integrands of $s'$ and $S''$,
but not as an explicit argument of $G$:
we now show that
it is redundant and can be eliminated.

From the Fundamental Equality for $f'$ and the pair $c_n, \hat C$
(note that the initial value $f'(x_0, c_1, \dots, c_n, \hat c, \hat C)$
does not depend on $c_n$)
we find
after division by factors $\partial_2{\it\Theta}$, $\partial_n G \neq 0$,
$$
  0
  = \partial_{1+n} s'(x, c_1, \dots, c_n)
  (S''(x, c_1, \dots, c_n, \hat c) + \hat C)
  + \partial_{1+n} S''(x, c_1, \dots, c_n, \hat c).
$$
Making use of the fact that this holds for all $\hat C\in \R$ and all
$x\in I$, we conclude

$0 = \partial_{1+n} \varphi'(x, c_1, \dots, c_n)$, and
$0 = \partial_{1+n} {\it\Phi}''(x, c_1, \dots, c_n, \hat c)$.

Therefore, setting
$\hat\varphi(x, c_1, \dots, c_{n-1}) := \varphi'(x, c_1, \dots, c_{n-1}, 0)$
and
${\it\hat\Phi}(x, c_1, \dots, c_{n-1}, \hat c)
 := {\it\Phi}''(x, c_1, \dots, c_{n-1}, 0, \hat c)$,
we obtain quadratures $\hat s$, $\hat S$ with the required properties.
\hfill $\square$

\medskip
Starting with an arbitrary effectively one-parametric integral in finite form,
we can apply Props. 3.4 and 3.5 in the following iterative algorithm
to eliminate all except at most two of the quadratures:

\medskip
\leftskip=.75in
\parindent=-.75in
{\bf Step 1.} \quad
{\bf if} there is only one quadrature, {\bf stop.}

\noindent
{\bf else} apply Prop. 3.4;

{\bf if} this results in case 1, {\bf repeat} step 1.

{\bf else} proceed with

{\bf Step 2.} \quad
{\bf if} there are only two quadratures, {\bf stop.}

\noindent
{\bf else} apply Prop. 3.5; {\bf repeat} step 2.

\leftskip=0in
\parindent=0in

\medskip
Depending on whether
this algorithm terminates in step 1 or 2
(which it must
since the number of quadratures is finite at the beginning and
decreases by 1 at each step
except in the one single step
in which Prop. 3.4 winds up in case 2),
one of the following situations is reached:

\smallskip
1. The original integral is equivalent to the integral by quadratures
$$
  \hat f(x, C) = {\it\Theta}(x, G(S(x) + C))
  \qquad (x\in I, C \in \R).
$$
This clearly is of the form stated in Thm 1.6,
with ${\it\hat\Theta}(x, c) := {\it\Theta}(x, G(c))$ $(x\in I, c\in \R)$,
$p := 0$, and $q := {\it\Phi}$ (the integrand of $S$).

\smallskip
2. The original integral is equivalent to the integral by quadratures
$$
  \hat f(x, c, C) = {\it\Theta}(x, G(e^{s(x) + c} (S(x, c) + C)))
  \qquad (x\in I, c, C \in \R).
$$
The Fundamental Equality for $\hat f$ and the pair $c, C$
yields after cancellation of factors $\partial_2{\it\Theta}$,
$\partial_1 G \neq 0$
$$
  0 = S(x, c) + \partial_2 S(x, c)
  \qquad\qquad (x\in I, c \in \R).
$$
Differentiating with respect to $x$ and using that $c$ is arbitrary,
we find
$$
  0 = {\it\Phi}(x, c) + \partial_2 {\it\Phi}(x, c)
  \qquad\qquad (x\in I, c \in \R).
$$
Thus there is a locally integrable function $q : I \rightarrow \R$
such that ${\it\Phi}(x, c) = q(x)\, e^{-c}$.
Eliminating the redundant integration constant of the quadrature $s$,
we arrive at the desired expression.

This completes the proof of Thm. 1.6.

\bigskip
{\bf Appendix. Integration by quadratures of the linear second-order equation.}

\bigskip
As is well known, the homogeneous linear second-order equation
$$
  u'' + Q\, u = 0,
\eqno (5)
$$
$Q : I \rightarrow \R$ locally integrable,
is equivalent to the first-order system
$$\eqalignno{
  \vartheta' &= \cos^2 \vartheta + Q\, \sin^2 \vartheta & (6)
\cr
  (\log\varrho)' &= (Q - 1)\, \sin\vartheta\,\cos\vartheta & (7)
\cr}$$
by the Pr\"ufer transformation
$$
  \pmatrix{u \cr u' \cr}
  = \varrho\, \pmatrix{\sin\vartheta \cr \cos\vartheta \cr}.
$$
Clearly the solution of $(7)$ can be obtained
from that of $(6)$ by a simple integration;
therefore it seems reasonable
to say that $(5)$ can be {\it integrated by quadratures\/} if and only if
$(6)$ can be integrated by quadratures in the sense specified in Section 1.

We call the family
$$
  f(x, c_1, \dots, c_n)
   := F(s_1(x) + c_1, \dots, s_n(x, c_1, \dots, c_{n-1}) + c_n)
$$
$(x\in I, c_1, \dots, c_n \in \R)$,
with $(s_1, \dots, s_n)$ an independent system of quadratures
and $F : \R^n \rightarrow \R$ admissible,
an {\it integral by quadratures in the restricted sense\/}
(excluding any explicit dependence of the integral
on the independent variable $x$).
This restriction narrows the concept of integration by quadratures considerably;
however, it is interesting to observe that then
the following consequence can be drawn from Thm. 1.6
without any assumptions of elementarity of the functions involved:

\medskip
{\bf 4.1 Theorem.}\qquad
{\it
If
the homogeneous linear second-order equation $(5)$
can be integrated by quadratures in the restricted sense,
then $Q$ is constant.
}

\medskip
{\it Proof.}\qquad
The right hand side of $(6)$ is uniformly Lipschitz with respect to $\vartheta$;
we can therefore apply Thm. 1.6
(with ${\it\Theta}(x, c) = c$) to learn that
$(6)$ can be integrated by quadratures in the restricted sense only if
there are a (four times continuously differentiable) diffeomorphism
$\phi : \R \rightarrow \R$,
and locally integrable functions $p, q$ such that
$$
  (\phi \circ \vartheta)'(x) + \phi \circ \vartheta(x)\, p(x) = q(x)
  \qquad (x\in I)
$$
whenever $\vartheta : I \rightarrow \R$ is a solution of $(6)$.
Applying the chain rule and noting that
through every point $(x, y) \in I \times \R$
there passes some solution of $(6)$, we conclude that
$$
  \phi'(y)\, (\cos^2 y + Q(x)\, \sin^2 y)
  + \phi(y)\, p(x) = q(x)
  \qquad (x\in I, y\in \R).
$$

Differentiating this identity with respect to y,
dividing by $\phi'(y) \neq 0$,
and differentiating again with respect to y
in order to eliminate $p$ and $q$, we find
$$\eqalignno{
  \Big({\phi'' \over \phi'}\Big)' (y)
  \, (\cos^2 y &+ Q(x)\, \sin^2 y)
  + {\phi'' \over \phi'}(y)
  \, ((\cos^2)'(y) + Q(x)\, (\sin^2)'(y)) &
\cr
  &+ (\cos^2)''(y) + Q(x)\, (\sin^2)''(y)
  = 0
  \qquad (x\in I, y \in \R). & (8)
\cr}$$

If $Q$ is not constant,
there are $x_1, x_2\in I$ such that $Q(x_1) \neq Q(x_2)$.
Taking the difference of $(8)$ at $x=x_1$ and $x=x_2$, it follows that
$$
  \Big({\phi'' \over \phi'}\Big)' (y) \,\sin^2 y
  + {\phi'' \over \phi'}(y) \, (\sin^2)'(y)
  + (\sin^2)''(y)
  = 0
  \qquad (y \in \R).
\eqno (9)
$$
When we multiply this equation by $1 - Q(x)$
and add the result to $(8)$,
we obtain
$$
  \Big({\phi'' \over \phi'}\Big)' = 0,
$$
which by $(9)$ implies that
$(\sin^2)'{} = 2\,\sin{}\cos$ and $(\sin^2)'' = 2\,(\cos^2 - \sin^2)$
are linearly dependent,
which is not the case.
\hfill $\square$

\bigskip
{\bf Biographical Note.}

\bigskip
The mathematical abilities of the young Vladimir Pavlovi\v c Maximovi\v c
early attracted \v Ceby\v sev's attention.
After leaving the 1st class at the physical-mathematical faculty of
St. Petersburg University in 1867, he completed his mathematical education as a
student of Bertrand at the Paris \'Ecole Polytechnique, where he obtained
in 1879 the doctoral degree for his thesis `Nouvelle m\'ethode pour int\'egrer
les equations simultan\'ees aux differentielles totales'.
His 1880 paper [3] was motivated by his hope to find, based on the principle
of the confluence, or mutual compensation, of integration constants, a proof
for the impossibility of integrating the general linear second-order equation
by quadratures.
These investigations, of which he had, already on 1st July, 1880, deposited a
preliminary version `M\'emoire sur les \'equations diff\'erentielles
g\'en\'erales du premier ordre qui s'int\`egrent au moyen d'un nombre fini de
quadratures. D\'emonstration de l'impossibilit\'e d'une telle int\'egration de
l'\'equation lin\'eaire du second ordre' at the Paris Academy in a sealed
envelope under the motto {\it nihil optimum nisi mathesis, et non est mortale
quod opto\/}, are published in his 1885 habilitation thesis [4] at Kazan',
where he was subsequently {\it Privatdozent\/}, and later {\it Dozent} at the
chair for pure mathematics at the physical-mathematical faculty.
There he published six more papers on the integration of differential equations,
on interpolation, function expansions and the roots of algebraic equations.
After he had moved to Kiev, Maximovi\v c's interest turned to probability
theory, and to the construction of a computing machine.
His last published work, and his only one on probability theory, was a talk
(Kiev 1888) on the application of probabilistic laws to school statistics,
using empirical data from the admission exams of the \'Ecole Polytechnique.
Early in the following year, signs of severe mental illness began to manifest
themselves, which led to the premature death in his 40th year, in St.
Petersburg on 17 October, 1889, of `the talented Russian mathematician whose
name in the history of science will be inseparably linked to the important
question of the integration of differential equations by quadratures'
([9] p.55 seq.).

\bigskip
{\it Acknowledgement.}\qquad
The author is indebted to H. Kalf for pointing out, and procuring copies of,
the rare publications [4] and [9], and for a thorough reading of the manuscript.
He also likes to thank J. Walter (Aachen) for his comments.

\bigskip
{\bf References.}

\bigskip
\leftskip=.25in
\parindent=-.25in
 
1. R Courant, D Hilbert: {\it Methoden der mathematischen Physik II.}
     Springer, Berlin 1937

2. S Lie: {\it Gesammelte Abhandlungen III.} Teubner, Leipzig 1922

3. W de Maximovitch: Conditions pour che les constantes arbitraires d'une
     expression g\'en\'erale soient distinctes entre elles.
     {\it Liouville J. (J. Math. Pures Appl.)} (3) {\bf VI} (1880) 167--177

4. V P Maximovi\v c: {\it Razyskanie ob\v s\v cyh differencial'nyh
     uravneni\u\i\ pervago poryadka, integriruyu\v s\v cyhsya v kone\v cnom
     vide, i dokazatel'stvo nevozmo\v znosti takogo integrirovaniya dlya
     ob\v s\v cago line\u\i nago uravneniya vtorago poryadka.}
     Tip. Imp. Univ., Kazan' 1885

5. P J Olver: {\it Applications of Lie groups to differential equations.}
     Springer, New York 1986

6. J F Pommaret: {\it Differential Galois theory.}
     Gordon \& Breach, New York 1983

7. J F Ritt: {\it Integration in finite terms.} Columbia Univ. Press,
     New York 1948

8. A V Vasil'ev. {\it Fortschr. d. Math.} {\bf 17} (1885) 305--309

9. A V Vasil'ev. {\it Sobranie protokolov zasedani\u\i\ sekci\u\i\ 
     fiziko-matemati\v ceskih nauk ob\v s\v cestva estestvoispytatele\u\i\ 
     pri Imperatorskom Kazanskom Universitete} {\bf 8} (1890) 53--56

\parindent=0in
\leftskip=0in

\bye